\documentclass[reqno,10pt, centertags, draft]{amsart}
\usepackage{amsmath,amsthm,amscd,amssymb,latexsym,upref}

\newcommand{\bbN}{{\mathbb{N}}}
\newcommand{\bbR}{{\mathbb{R}}}

\newcommand{\bbC}{{\mathbb{C}}}

\newcommand{\cC}{{\mathcal C}}

\newcommand{\cH}{{\mathcal H}}

\newcommand{\cJ}{{\mathcal J}}

\newcommand{\no}{\nonumber}
\newcommand{\lb}{\label}

\newcommand{\ol}{\overline}

\newcommand{\wti}{\widetilde}

\newcommand{\loc}{\text{\rm{loc}}}

\newcommand{\dom}{\text{\rm{dom}}}

\newcommand{\supp}{\text{\rm{supp}}}

\newcommand{\bi}{\bibitem}

\renewcommand{\bar}[1]{\overline{#1}}

\newtheorem{theorem}{Theorem}

\newtheorem{hypothesis}[theorem]{Hypothesis}
\theoremstyle{definition}

\begin{document}

\title[$\cJ$-Self-Adjointness of a Class of Dirac-type Operators]
{$\cJ$-Self-Adjointness of a Class of Dirac-type Operators}
\author[R.\ Cascaval and F.\ Gesztesy]{Radu
Cascaval and Fritz Gesztesy}
\address{Department of Mathematics,
University of Colorado, Colorado Springs, CO 80933}
\email{radu@math.uccs.edu}
\urladdr{http://math.uccs.edu/$\sim $radu}
\address{Department of Mathematics,
University of Missouri, Columbia, MO 65211, USA}
\email{fritz@math.missouri.edu}
\urladdr{http://www.math.missouri.edu/people/fgesztesy.html}
\thanks{To appear in {\it J. Math. Anal. Appl.}}


\begin{abstract}
In this note we prove that the maximally defined operator associated with
the Dirac-type differential expression
$M(Q)=i\left(\begin{matrix} \frac{d}{dx}I_m &
-Q\\ -Q^* & -\frac{d}{dx}I_m\end{matrix}\right)$, where $Q$ represents
a symmetric $m\times m$ matrix (i.e., $Q(x)^\top=Q(x)$ a.e.) with entries
in $L^1_{\loc}(\mathbb{R})$, is $\mathcal{J}$-self-adjoint,
where $\mathcal{J}$ is the antilinear conjugation defined by
$\mathcal{J}=\sigma_1\mathcal{C}$, $\sigma_1= \left(\begin{smallmatrix}
0 & I_m\\ I_m & 0 \end{smallmatrix}\right)$ and
$\mathcal{C}(a_1,\ldots,a_m, b_1,\ldots, b_m)^{\top}=
(\bar{a_1},\ldots,\bar{a_m}, \bar{b_1},\ldots, \bar{b_m})^{\top}$.
The differential expression $M(Q)$ is of significance as it appears in the
Lax formulation of the nonabelian (matrix-valued) focusing nonlinear
Schr\"odinger hierarchy of evolution equations.
\end{abstract}

\maketitle

To set the stage for this note, we briefly mention the Lax pair and
zero-curvature representations of the matrix-valued
Ablowitz-Kaup-Newell-Segur (AKNS) equations and the special focusing and
defocusing nonlinear Schr\"odinger (NLS) equations associated with it. Let
$P=P(x,t)$ and $Q=Q(x,t)$ be smooth $m\times m$ matrices, $m\in\bbN$, and
introduce the Lax pair of $2m\times 2m$ matrix-valued differential
expressions
\begin{align}
M(P,Q)&=i\begin{pmatrix} \frac{d}{dx}I_m &
-Q\\ P & -\frac{d}{dx}I_m \end{pmatrix} \lb{0.1} \\
L(P,Q)&=i\begin{pmatrix} \frac{d^2}{dx^2}I_m-\frac{1}{2}QP &
-Q\frac{d}{dx}-\frac{1}{2}Q_x\\ P\frac{d}{dx}-\frac{1}{2}P_x &
-\frac{d^2}{dx^2}I_m+\frac{1}{2}PQ\end{pmatrix} \lb{0.2}
\end{align}
and the $2m\times 2m$ zero-curvature matrices
\begin{align}
U(z,P,Q)&=\begin{pmatrix} -izI_m & Q \\ P & izI_m\end{pmatrix}, \lb{0.3}
\\
V(z,P,Q)&=\begin{pmatrix} -iz^2I_m-\frac{i}{2}QP &
zQ+\frac{i}{2}Q_x \\ zP-\frac{i}{2}P_x &
iz^2I_m+\frac{i}{2}PQ\end{pmatrix}, \lb{0.4}
\end{align}
where $z\in\bbC$ denotes a (spectral) parameter and $I_m$ is the
identity matrix in $\bbC^m$. Then the Lax equation
\begin{equation}
\frac{d}{dt}M-[L,M]=0 \lb{0.5}
\end{equation}
is equivalent to the $m\times m$ matrix-valued AKNS system
\begin{align}
\begin{split}
Q_t-\frac{i}{2}Q_{xx}+iQPQ&=0, \\
P_t+\frac{i}{2}P_{xx}-iPQP&=0, \lb{0.6}
\end{split}
\end{align}
where $[\cdot,\cdot]$ denotes the commutator symbol. Similarly, the
zero-curvature equation
\begin{equation}
U_t-V_x+[U,V]=0  \lb{0.7}
\end{equation}
is also equivalent to the $m\times m$ matrix-valued AKNS system
\eqref{0.6}. Two special cases of this formalism are of particular
importance: The  {\it focusing} NLS equation
\begin{align}
\textit{focusing: } \;\; &Q_t-\frac{i}{2}Q_{xx}-iQQ^*Q=0,  \lb{0.8} \\
\intertext{obtained from \eqref{0.1}--\eqref{0.7} in the special case
where $P=-Q^*$, and the {\it defocusing} NLS equation}
\textit{defocusing: } \;\; &Q_t-\frac{i}{2}Q_{xx}+iQQ^*Q=0,  \lb{0.9}
\end{align}
obtained from \eqref{0.1}--\eqref{0.7} in the special case where
$P=Q^*$. Here $Q^*$ denotes the adjoint (i.e., complex conjugate and
transpose) matrix of $Q$.

In this note we will restrict our attention to the focusing NLS case
$P=-Q^*$. (See, e.g., \cite[Sect.\ 3.3, Ch.\ 8]{AK90}, \cite{TW98} and
\cite[Sect 3.1]{Ts00} in which an inverse  scattering approach is
developed for the matrix NLS equation \eqref{0.8}). Actually, \eqref{0.6},
\eqref{0.8}, and \eqref{0.9} are just the first equations in an infinite
hierarchy of nonlinear evolution equations (the nonabelian AKNS, and
focusing and defocusing NLS hierarchies) but we will not further dwell on
this point.

Recently, it has been proved in \cite{CGHL03}  that the $2\times 2$
matrix-valued Lax differential  expression
\begin{equation}
M(q)=i\left(\begin{matrix} \frac{d}{dx} &
-q\\ -\bar{q} & -\frac{d}{dx}\end{matrix}\right)\lb{1.1}
\end{equation}
corresponding to the scalar focusing NLS hierarchy defines, under the most
general hypothesis $q\in L^1_{\loc}(\bbR)$ on the potential $q$, a
$\cJ$-self-adjoint operator in $L^2(\bbR)^2$, where $J$ is the antilinear
conjugation, $\cJ=\sigma_1 \cC$, with $\sigma_1=\left(\begin{smallmatrix}
0 & 1\\ 1 &0 \end{smallmatrix}\right)$ and $\cC$ is complex conjugation in
$\bbC^2$.  This is the direct analog of a recently proven fact in
\cite[Lemma\ 2.15]{CG02} that the Dirac-type  Lax differential expression
in the defocusing nonlinear Schr\"odinger (NLS) case is always in the
limit point case at $\pm\infty$. Equivalently, the maximally defined
Dirac-type operator corresponding to the defocusing NLS case is always
self-adjoint.

In this paper we present an extension of the result in \cite{CGHL03}, for
$2m\times 2m$ matrix-valued Dirac-type differential expressions of the
form
\begin{equation}
M(Q)=i\left(\begin{matrix} \frac{d}{dx}I_m &
-Q\\ -Q^* & -\frac{d}{dx}I_m \end{matrix}\right) \lb{1.2}
\end{equation}
associated with the nonabelian (matrix-valued) focusing NLS equation
\eqref{0.8}.

We will assume the following conditions on $Q$ from now on ($A^\top$
denotes the transpose of the matrix $A$):

\begin{hypothesis}\label{h.1}
Assume $Q\in L^1_\loc(\bbR)^{m\times m}$ satisfies
\begin{equation}
	Q=Q^\top \text{ a.e.} \lb{1.5}
\end{equation}
\end{hypothesis}

Next, we briefly recall some basic facts about $\cJ $-symmetric
and $\cJ $-self-adjoint operators in a complex Hilbert space $\cH$ (see, e.g.,
\cite[Sect.\ III.5]{EE89} and \cite[p.\ 76]{Gl65}) with
scalar product denoted by $(\cdot ,\cdot)_\cH$ (linear in the first
and antilinear in the second place) and corresponding norm denoted by
$\|\cdot\|_\cH$. Let $\cJ$ be a conjugation operator in $\cH$, that is,
$\cJ $ is an antilinear involution satisfying
\begin{equation}
	(\cJ u,v)_\cH=(\cJ v,u)_\cH \text{ for all }u,v\in \cH,\quad \cJ^2=I.
\lb{1.6}
\end{equation}
In particular,
\begin{equation}
	(\cJ u,\cJ v)_\cH=(v,u)_\cH, \quad u,v\in \cH. \lb{1.7}
\end{equation}
A linear operator $S$ in $\cH$, with domain $\dom(S)$ dense in $\cH$,
is called $\cJ $-symmetric if
\begin{equation}
S\subseteq \cJ  S^*\cJ \, \text{ (equivalently, if
$JSJ\subseteq S^*$).}  \lb{1.8}
\end{equation}
Clearly, \eqref{1.8} is equivalent to
\begin{equation}
	(\cJ u,Sv)_\cH=(\cJ Su,v)_\cH, \quad u,v\in \dom(S). \lb{1.9}
\end{equation}
Here $S^*$ denotes the adjoint operator of $S$ in $\cH$. If $S$ is
$\cJ$-symmetric, so is its closure $\ol{S}$. The operator
$S$ is called $\cJ $-self-adjoint if
\begin{equation}
	S=\cJ S^*\cJ \, \text{ (equivalently, if $\cJ S\cJ =S^*$).}
\end{equation}
Finally, a densely defined, closable operator $T$ is called essentially
$\cJ$-self-adjoint if its closure, $\ol{T}$, is $\cJ $-self-adjoint,
that is, if
\begin{equation}
\ol{T}=\cJ T^* \cJ.
\end{equation}
Next, assuming $S$ to be $\cJ$-symmetric, one introduces the following
inner product $(\cdot,\cdot)_*$ on $\dom(\cJ S^*\cJ )=\cJ\dom(S^*)$
according to
\cite{Kn81} (see also \cite{Ra85}),
\begin{equation}
	(u,v)_*=(\cJ u,\cJ v)_\cH+(S^*\cJ u,S^*\cJ v)_\cH,
\quad u,v\in \dom(\cJ S^*\cJ),
\end{equation}
which renders $\dom(\cJ S^*\cJ )$ a Hilbert space. Then the following
theorem holds ($I_\cH$ denotes the identity operator in $\cH$).
\begin{theorem}[Race  \cite{Ra85}] \lb{t.2}
Let $S$ be a densely defined closed $\cJ $-symmetric operator. Then
\begin{equation}
	\dom (\cJ S^*\cJ )=\dom (S) \oplus_* \ker((S^*\cJ )^2+I_\cH), \lb{1.10}
\end{equation}
where $\oplus_*$ means the orthogonal direct sum with respect to the inner
product $(\cdot,\cdot)_*$. In particular, a densely defined closed $\cJ
$-symmetric operator $S$ is $\cJ $-self-adjoint if and only if
\begin{equation}
	\ker((S^*\cJ )^2+I_\cH)=\{ 0\} .\lb{1.11}
\end{equation}
\end{theorem}

Theorem \ref{t.2} will be used to prove the principal result of this note
that the (maximally defined) Dirac-type operator associated with the
differential expression $M(Q)$ in \eqref{1.2} (relevant to the focusing
matrix NLS equation \eqref{0.8}) is  always $\cJ$-self-adjoint under most
general conditions on the coefficient $Q$ in Hypothesis \ref{h.1} (see
Theorem \ref{t.5} below). This will be done by verifying a  relation of
the type \eqref{1.11}.

To this end, it is convenient to introduce some
standard notations to be used throughout
the remainder of this paper. The Hilbert space $\cH$ is chosen to be
$L^2 (\bbR)^{2m}=L^2 (\bbR)^{m}\oplus L^2 (\bbR)^{m}$.
The space of $m\times m$ matrices with  entries in $L^1_\loc(\bbR)$
is denoted by $L^1_\loc(\bbR)^{m\times m}$. An antilinear conjugation
$\cJ$ in the complex Hilbert space $L^2(\bbR)^{2m}$ is defined by
\begin{equation}
\mathcal{J}=\sigma_1\mathcal{C}, \lb{J}
\end{equation}
where
\begin{equation}
  \sigma_1= \begin{pmatrix}
0 & I_m\\ I_m & 0 \end{pmatrix}, \quad
\mathcal{C}(a_1,\ldots,a_m, b_1,\ldots, b_m)^{\top}=
(\bar{a_1},\ldots,\bar{a_m}, \bar{b_1},\ldots, \bar{b_m})^{\top}.
\end{equation}

Given Hypothesis \ref{h.1}, we now introduce the following maximal and
minimal Dirac-type operators in $L^2(\bbR)^{2m}$ associated with the
differential expression $M(Q)$,
\begin{align}
&D_{\max}(Q)F=M(Q)F, \lb{1.12} \\
&F\in\dom(D_{\max}(Q))=\big\{G\in L^2(\bbR)^{2m}\,\big|\, G\in
AC_\loc(\bbR)^{2m}, \, M(Q)G\in L^2(\bbR)^{2m}\big\}, \no \\
&D_{\min}(Q)F=M(Q)F,  \lb{1.13} \\
&F\in\dom(D_{\min}(Q))=\{G\in\dom(D_{\max}(Q)) \,|\,
\supp(G) \text{ is compact}\}. \no
\end{align}
It follows by standard techniques (see, e.g., \cite[Ch.\ 8]{LS91} and
\cite{We71}) that under Hypothesis~\ref{h.1},
$D_{\min}(Q)$ is densely defined and closable in $L^2(\bbR)^{2m}$ and
$D_{\max}(Q)$ is a densely defined closed operator in $L^2(\bbR)^{2m}$.
Moreover one infers (see, e.g., \cite[Lemma
8.6.2]{LS91} and \cite{We71} in the analogous case of symmetric Dirac
operators)
\begin{equation}
\ol{D_{\min}(Q)}=D_{\max}(-Q)^*,  \, \text{ or equivalently, }\,
D_{\min}(Q)^*=D_{\max}(-Q). \lb{1.14}
\end{equation}

The following result will be a crucial ingredient in the proof of
Theorem \ref{t.5}, the principal result of this note.

\begin{theorem}\lb{t.3}
Assume Hypothesis \ref{h.1}.
Let $N(Q)$ be the following $($\,formally self-adjoint\,$)$
differential expression
\begin{equation}
	N(Q)=i\begin{pmatrix} \frac{d}{dx}I_m &
	-Q\\Q^*&\frac{d}{dx}I_m\end{pmatrix} \lb{1.15}
\end{equation}
and denote by $\wti D_{\max}(Q)$ the maximally defined Dirac-type operator
in $L^2(\bbR)^{2m}$ associated with $N(Q)$,
\begin{align}
&\wti D_{\max}(Q)F=N(Q)F, \lb{1.16} \\
&F\in\dom(\wti D_{\max}(Q))=\big\{G\in L^2(\bbR)^{2m}\,\big|\, G\in
AC_\loc(\bbR)^{2m}, \, N(Q)G\in L^2(\bbR)^{2m}\big\}. \no
\end{align}
Then, \\
$(i)$ The following identity holds
\begin{equation}
	M(-Q)M(Q)=N(Q)^2.\lb{1.17}
\end{equation}
$(ii)$ Let $U_Q=U_Q(x)$ satisfy the initial value problem
\begin{equation}
	U_Q'= \begin{pmatrix} 0 & Q\\ -Q^*& 0 \end{pmatrix} U_Q,
\quad U_Q(0)=I_{2m}. \lb{1.18}
\end{equation}
Then $\{U_Q(x)\}_{x\in \bbR}$ is a family of unitary matrices in
$\bbC^{2m}$ with entries in $AC_{\loc}(\bbR)\cap L^\infty(\bbR)$
satisfying
\begin{equation}
	U_Q^{-1}N(Q)U_Q=i\frac{d}{dx}I_{2m}.\lb{1.19}
\end{equation}
$(iii)$ Let $\mathcal{U}_Q$ denote the multiplication operator with $
U_Q(\cdot )$ on $L^2(\bbR)^{2m}$. Then $\wti D_{\max}(Q)$ is unitarily
equivalent to the maximally defined operator in $L^2(\bbR)^{2m}$
associated with the differential expression $i\frac{d}{dx}I_{2m}$,
\begin{align}
	&\mathcal{U}_Q^{-1}\wti D_{\max}(Q)\mathcal{U}_Q=\bigg(
	i\frac{d}{dx}I_{2m}\bigg)_{\max},\lb{1.20}\\
&	\dom \bigg(\bigg( i\frac{d}{dx}I_{2m}\bigg)_{\max}\bigg)=
H^{1,2}(\bbR)^{2m}  \no\\
& \quad =\big\{ F\in L^2(\bbR)^{2m}\, \big| \, F\in AC_{\loc}(\bbR)^{2m},
\, F'\in L^2(\bbR)^{2m}\big\}. \no
\end{align}
Moreover,
\begin{align}
&\mathcal{U}_Q^{-1}D_{\max}(-Q)D_{\max}(Q)\mathcal{U}_Q =
	\bigg( -\frac{d^2}{dx^2}I_{2m}\bigg)_{\max}, \lb{1.21} \\
&\dom \bigg(\bigg( -\frac{d^2}{dx^2}I_{2m}\bigg)_{\max}\bigg)
= H^{2,2}(\bbR)^{2m}\no\\
&\quad =\big\{F\in L^2(\bbR)^{2m}\, \big|\, F,F'\in AC_{\loc}(\bbR)^{2m},
\, F', F''\in L^2(\bbR)^{2m}\big\}. \no
\end{align}
\end{theorem}
\begin{proof}
That $N(Q)$ is formally self-adjoint and $M(-Q)M(Q)=N(Q)^2$, as stated in
$(i)$, is an elementary matrix calculation. \\
To prove $(ii)$, we note that the initial value problem
\eqref{1.18} is well-posed in the sense of Carath\'eodory since $Q\in
L^1_{\loc}(\bbR)^{m\times m}$
(cf., e.g., \cite[Lemma IX.2.2]{GGK90})  with a solution
matrix $U_Q$ with entries in $AC_\loc(\bbR)$. Moreover, for each
$x\in\bbR$, $U_Q(x)$ is a unitary matrix in $\bbC^{2m}$, since
$U_Q'=-B(Q)U_Q$, with
$B(Q)=\left(\begin{smallmatrix} 0 & -Q \\ Q^* & 0
\end{smallmatrix}\right)$ being skew-adjoint. Thus, the entries
$U_{Q,j,k}$, $1\leq j,k \leq {2m}$ of $U_Q$ (as well as those of
$U_Q^{-1}$) actually satisfy
\begin{equation}
U_{Q,j,k}\in AC_{\loc}(\bbR)\cap L^\infty(\bbR), \quad 1 \leq j,k \leq 2m.
\lb{1.22}
\end{equation}
(Since $U_Q$ is a bounded matrix-valued operator of multiplication in
$L^2(\bbR)^{2m}$, its entries $U_{Q,j.k}$ are all in $L^\infty(\bbR)$,
as one readily verifies by studying scalar products of the form
$(F_j,U_QF_k)_{L^2(\bbR)^{2m}}=(f_j,U_{Q,j,k}f_k)_{L^2(\bbR)}$,
$1\leq j,k\leq 2m$, where $F_j=(0,\dots,0,f_j,0,\dots,0)^\top$ with
$f_j\in L^2(\bbR)$, $1\leq j\leq 2m$.) Next, fix
$F\in AC_{\loc}(\bbR)^{2m}$, such that $\mathcal{U}_Q^{-1}F\in
H^{1,2}(\bbR)^{2m}$. Then
\begin{align}
U_Q\bigg( i\frac{d}{dx}I_{2m}\bigg) U_Q^{-1}F&=i\frac{d}{dx}F+ iU_Q
\frac{d}{dx}(U_Q^{-1})F \no\\
&=i\frac{d}{dx}F+iU_Q(U_Q^{-1}B(Q)^*)F \no\\
&=N(Q)F, \lb{1.23}
\end{align}
where we used the fact that  $(U_Q^{-1})'=U_Q^{-1}B(Q)^*$.
Thus, $(ii)$ follows. \\
Moreover, by \eqref{1.22}, the fact that $U_Q$ is unitary in
$\bbC^{2m}$, and by \eqref{1.23} one concludes
$\dom(\wti D_{\max}(Q))= \mathcal{U}_Q H^{1,2}(\bbR)^{2m}$ . This proves
\eqref{1.20}. \\
Clearly $(i)$ and $(ii)$ yield the relation
\begin{equation}
U_Q^{-1}M(-Q)M(Q)U_Q = -\frac{d^2}{dx^2}I_{2m}. \no
\end{equation}
Thus, \eqref{1.21} will follow once we prove the following facts:
\begin{align}
& (i) \; U_Q F\in L^2(\bbR)^{2m} \text{ if and only if } F\in L^2(\bbR)^{2m},
\lb{1.24a} \\
& (ii) \; U_Q F\in AC_\loc (\bbR)^{2m} \text{ if and only if } F\in
AC_\loc(\bbR)^{2m}, \lb{1.24b} \\
& (iii) \; M(Q)U_Q F\in L^2(\bbR)^{2m} \text{ if and only if }
F'\in L^2(\bbR)^{2m}, \lb{1.24c} \\
& (iv) \; M(Q)U_Q F\in AC_\loc (\bbR)^{2m} \text{ if and only if }
F'\in AC_\loc(\bbR)^{2m}, \lb{1.24d} \\
& (v) \; M(-Q)M(Q)U_Q F\in L^2(\bbR)^{2m} \text{ if and only if }
F''\in L^2(\bbR)^{2m}. \lb{1.24e}
\end{align}
Clearly \eqref{1.24a} and \eqref{1.24e} hold since $U_Q$ is unitary in
$\bbC^{2m}$. \eqref{1.24b} is valid since \\
$U_{Q,j,k}, U^{-1}_{Q,j,k}\in
AC_{\loc}(\bbR)\cap L^\infty(\bbR)$, $j,k=1,\ldots ,2m$. Next, for
$F=(F_1^{\top}, F_2^{\top})^{\top}$, $F_1, F_2\in L^2(\bbR)^m$, an
explicit computation yields
\begin{equation}
M(Q)U_QF=i\begin{pmatrix} U_Q^{(1)} F_1'+U_Q^{(2)}F_2' \\
-U_Q^{(3)}F_1'-U_Q^{(4)}F_2' \end{pmatrix}, \quad F=(F_1^{\top},\;
  F_2^{\top})^\top, \lb{1.25}
\end{equation}
where $U_Q^{(l)}$, $l=1,2,3,4$ are blocks of the matrix $U_Q$,
\begin{equation}
U_Q=\begin{pmatrix} U_Q^{(1)}& U_Q^{(2)}\\ U_Q^{(3)}& U_Q^{(4)}
\end{pmatrix}.
\end{equation}
Introducing
\begin{equation}
V_Q=\sigma_3 U_Q\sigma_3
=\begin{pmatrix} U_Q^{(1)} & -U_Q^{(2)} \\ -U_Q^{(3)} & U_Q^{(4)}
\end{pmatrix}, \lb{1.26}
\end{equation}
one infers $V_{Q,j,k}, V^{-1}_{Q,j,k}\in
AC_{\loc}(\bbR)\cap L^\infty(\bbR)$, $j,k=1,\ldots, 2m$ and
\begin{equation}
V_Q^{-1}M(Q)U_QF=i(F_1', -F_2')^\top, \lb{1.27}
\end{equation}
and hence \eqref{1.24c} and \eqref{1.24d} hold. This proves \eqref{1.21}.
\end{proof}

The principal result of this note then reads as follows.

\begin{theorem} \lb{t.5}
Assume Hypothesis \ref{h.1}. Then the  minimally defined Dirac-type
operator $D_{\min}(Q)$ associated with the Lax differential expression
\begin{equation}
	M(Q)=i\begin{pmatrix} \frac{d}{dx}I_m &
	-Q\\-Q^*&-\frac{d}{dx}I_m\end{pmatrix}
\end{equation}
introduced in \eqref{1.13} is essentially $\cJ$-self-adjoint in
$L^2(\bbR)^{2m}$, that is,
\begin{equation}
\ol{D_{\min}(Q)}=\cJ D_{\min}(Q)^* \cJ, \lb{1.28}
\end{equation}
where $\cJ$ is the conjugation defined in \eqref{J}. Moreover,
\begin{equation}
\ol{D_{\min}(Q)}=D_{\max}(Q) \lb{1.29}
\end{equation}
and hence $D_{\max}(Q)$ is $\cJ$-self-adjoint.
\end{theorem}
\begin{proof}
We first recall (cf.\ \eqref{1.14})
\begin{equation}
D_{\min}(Q)^*=D_{\max}(-Q) \lb{1.30}
\end{equation}
and also note
\begin{equation}
\cJ D_{\max}(-Q)\cJ=D_{\max}(Q^\top)=D_{\max}(Q). \lb{1.31}
\end{equation}
Here we employed the symmetry of $Q$ (see \eqref{1.5}).
Since $\ol{D_{\min}(Q)}$ is closed and $\cJ$-symmetric
(this follows from \eqref{1.30} and \eqref{1.31}), its
$\cJ$-self-adjointness is equivalent to showing that
(cf. \eqref{1.11})
\begin{align}
&\ker \big(D_{\min} (Q)^* \cJ D_{\min}(Q)^* \cJ+I_{L^2(\bbR)^{2m}}\big)
\no \\
& \quad = 
	\ker\big(D_{\max}(-Q)D_{\max}(Q)+I_{L^2(\bbR)^{2m}}\big) 
=\{0\}.
\lb{1.32}
\end{align}

Since by Theorem \ref{t.3}\,(iii), $D_{\max}(-Q)D_{\max}(Q)$ is
unitarily equivalent to \\
$(-d^2/dx^2 I_{2m})_{\max}\geq 0$, one concludes that
\begin{equation}
D_{\max}(-Q)D_{\max}(Q)\geq 0\no
\end{equation}
and hence \eqref{1.32} holds.
The fact \eqref{1.29} now follows from \eqref{1.28} and
\eqref{1.30},
\begin{equation}
	\ol{D_{\min}(Q)}=\cJ D_{\min}(Q)^*\cJ =\cJ D_{\max}(-Q)\cJ
	=D_{\max}(Q).
\end{equation}
\end{proof}

As mentioned in the introductory paragraph, Theorem \ref{t.5} in the
$\cJ$-self-adjoint context can be viewed as an analog of Lemma\ 2.15 in
\cite{CG02} in connection with self-adjoint Dirac-type operator relevant
in the nonabelian (matrix-valued) defocusing nonlinear Schr\"odinger
hierarchy (cf.\ also \cite{LM02} for results of this type).

We conclude with a short remark. The special case where
\begin{equation}
Q=\begin{pmatrix} q_1 & 0 &\dots & 0\\
q_2 & 0 & \dots & 0 \\ \vdots & \vdots & \ddots &\vdots \\
q_m & 0 & \dots & 0 \end{pmatrix}, \, \text{ or } \;
Q=({\bf q \quad 0}),  \lb{1.3a}
\end{equation}
is known as the vector NLS equation (cf.\ \cite{Ch96})
\begin{equation}
	i{\bf q}_t+\frac12 {\bf q}_{xx}+\| {\bf q}\|^2 {\bf q}=0, \lb{1.4}
\end{equation}
a generalization of the well-known Manakov system \cite{Ma74} (for $m=2$).
Here ${\bf q}=(q_1, \ldots q_m)^t$, ($\| {\bf q}\|^2={\bf q}^*{\bf q}
=\sum_{j=1}^m |q_j|^2$). Unfortunately, the methods applied in this note
forced us to restrict our attention to symmetric matrices $Q$ only
(i.e., $Q=Q^\top$) and hence our current result does not apply to the vector
NLS case.


\end{document}